\documentclass[conference, 9pt]{IEEEtran}
\IEEEoverridecommandlockouts
\usepackage{graphicx}
\usepackage{color}
\usepackage{placeins}
\usepackage{float}
\usepackage{tabularx,colortbl}
\usepackage{amsmath, amssymb}
\usepackage{epsfig}

\usepackage{amsthm}

\begin{document}
\title{Conservative Signal Processing Architectures\\For Asynchronous, Distributed Optimization\\Part II: Example Systems}

\author{\IEEEauthorblockN{Thomas A. Baran and Tarek A. Lahlou\thanks{The authors wish to thank Analog Devices, Bose Corporation, and Texas Instruments for their support of innovative research at MIT and within the Digital Signal Processing Group.}}
\IEEEauthorblockA{Digital Signal Processing Group\\
Massachusetts Institute of Technology}}

\maketitle

\begin{abstract}

This paper provides examples of various synchronous and asynchronous signal processing systems for performing optimization, utilizing the framework and elements developed in a preceding paper.  The general strategy in that paper was to perform a linear transformation of stationarity conditions applicable to a class of convex and nonconvex optimization problems, resulting in algorithms that operate on a linear superposition of the associated primal and dual decision variables.  The examples in this paper address various specific optimization problems including the LASSO problem, minimax-optimal filter design, the decentralized training of a support vector machine classifier, and sparse filter design for acoustic equalization.  Where appropriate, multiple algorithms for solving the same optimization problem are presented, illustrating the use of the underlying framework in designing a variety of distinct classes of algorithms.  The examples are accompanied by numerical simulation and a discussion of convergence.
\end{abstract}

\IEEEpeerreviewmaketitle

\vspace{-2mm}
\section{Introduction}

This paper presents various classes of asynchronous, distributed optimization systems, demonstrating the use of the framework discussed in Part I \cite{BaranLahlouPartI}.  The design and use of each class of systems is based upon the following strategy:
\begin{enumerate}
\item Write a reduced-form optimization problem, defined in \cite{BaranLahlouPartI}.
\item Connect appropriate constitutive relations to interconnection elements, e.g.~from Figs.~2-3 in \cite{BaranLahlouPartI}, implementing the associated transformed stationarity conditions.  Delay-free loops will generally result.
\item Break delay-free loops:
  \begin{enumerate}
  \item For any constitutive relation that is a source element, perform algebraic simplification thereby incorporating the solution of the algebraic loop into the interconnection.
  \item Insert synchronous or asynchronous delays between the remaining constitutive relations and the interconnection.
  \end{enumerate}
\item Run the distributed system until it reaches a fixed point.  The discussion in Section \ref{sec:conv}, in conjunction with the system properties in Fig.~3 in \cite{BaranLahlouPartI}, provide guidance in determining when convergence is ensured.
\item Read out the primal and dual decision variables $a_i$ and $b_i$ by multiplying the variables $c_i$ and $d_i$ by the inverses of the ($2\times 2$) matrices used in transforming the stationarity conditions.
\end{enumerate}

\section{Example Systems}
\label{sec:examples}
Figs.~\ref{fig:lasso_approx}-\ref{fig:sparse_equalizer} depict various asynchronous, distributed optimization algorithms implemented using the presented framework, specifically making use of the elements in Figs.~2-3 of Part I \cite{BaranLahlouPartI}.
Figs.~\ref{fig:lasso_approx} and \ref{fig:lasso_augmented} in this paper illustrate two alternative implementations of systems for solving the LASSO problem. Figs.~\ref{fig:minimax_filter_design} and \ref{fig:minimax_filter_design_2} depict two alternative implementations of systems for performing minimax-optimal FIR filter design.  Fig.~\ref{fig:svm} depicts a support vector machine classifier trained using a decentralized algorithm generated using the presented framework.  Fig.~\ref{fig:sparse_equalizer} illustrates an example of a nonconvex optimization algorithm aimed at the problem discussed in \cite{baran2010linear}, in particular that of designing a sparse FIR filter for acoustic equalization.  In Figs.~\ref{fig:lasso_approx}-\ref{fig:sparse_equalizer}, the asynchronous delay elements were numerically simulated using discrete-time sample-and-hold systems triggered by independent Bernoulli processes, with the probability of sampling being $0.1$.

\section{Discussion of convergence}
\label{sec:conv}

\begin{figure}[h!]
\begin{center}
\epsfig{file=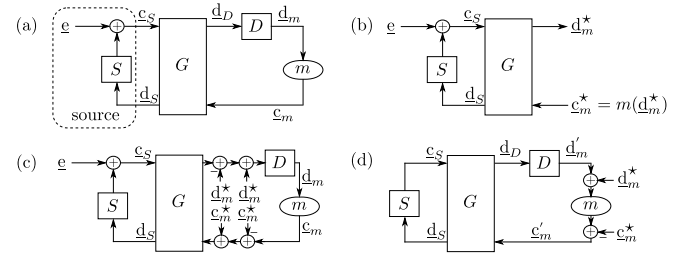,width=3.4in}
\vspace{-1mm}
\caption{(a) General description of the interconnection of elements used in the presented systems.  (b)-(d)  Manipulations performed in analyzing convergence.\label{fig:convergence}}
\end{center}
\end{figure}

\begin{figure*}[t!]
\begin{center}
\vspace{-4mm}
\epsfig{file=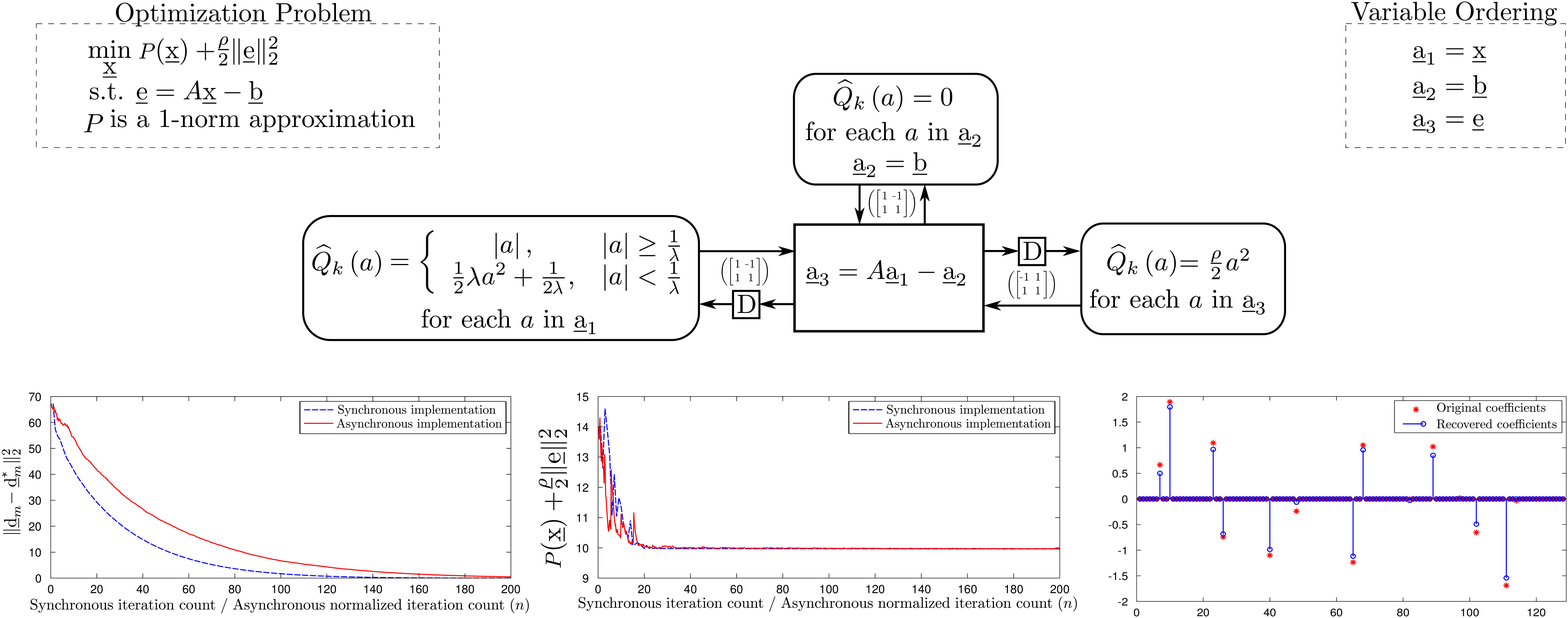,width=7in} 
\caption{Signal processing architecture and numerical simulation corresponding to an algorithm for solving the LASSO problem.  An approximation to the $1$-norm is used that is quadratic in the close vicinity of $0$.  The parameters $\lambda$ and  $\rho$ are selected to specify the interval outside of which the $1$-norm approximation is exact and to trade off between the sparsity of the solution and the enforcement of the linear equality constraints, respectively. For the depicted solution $\lambda$ and $\rho$ are selected to be large.  Note in particular the monotonic convergence of $|| \underbar{d}_m - \underbar{d}_m^{\star}||^2_2 $ to zero.  ``Asynchronous normalized iteration count'' indicates the number of iterations times the probability of sampling, discussed in Section \ref{sec:examples}.}\label{fig:lasso_approx}
\vspace{0.15in}
\epsfig{file=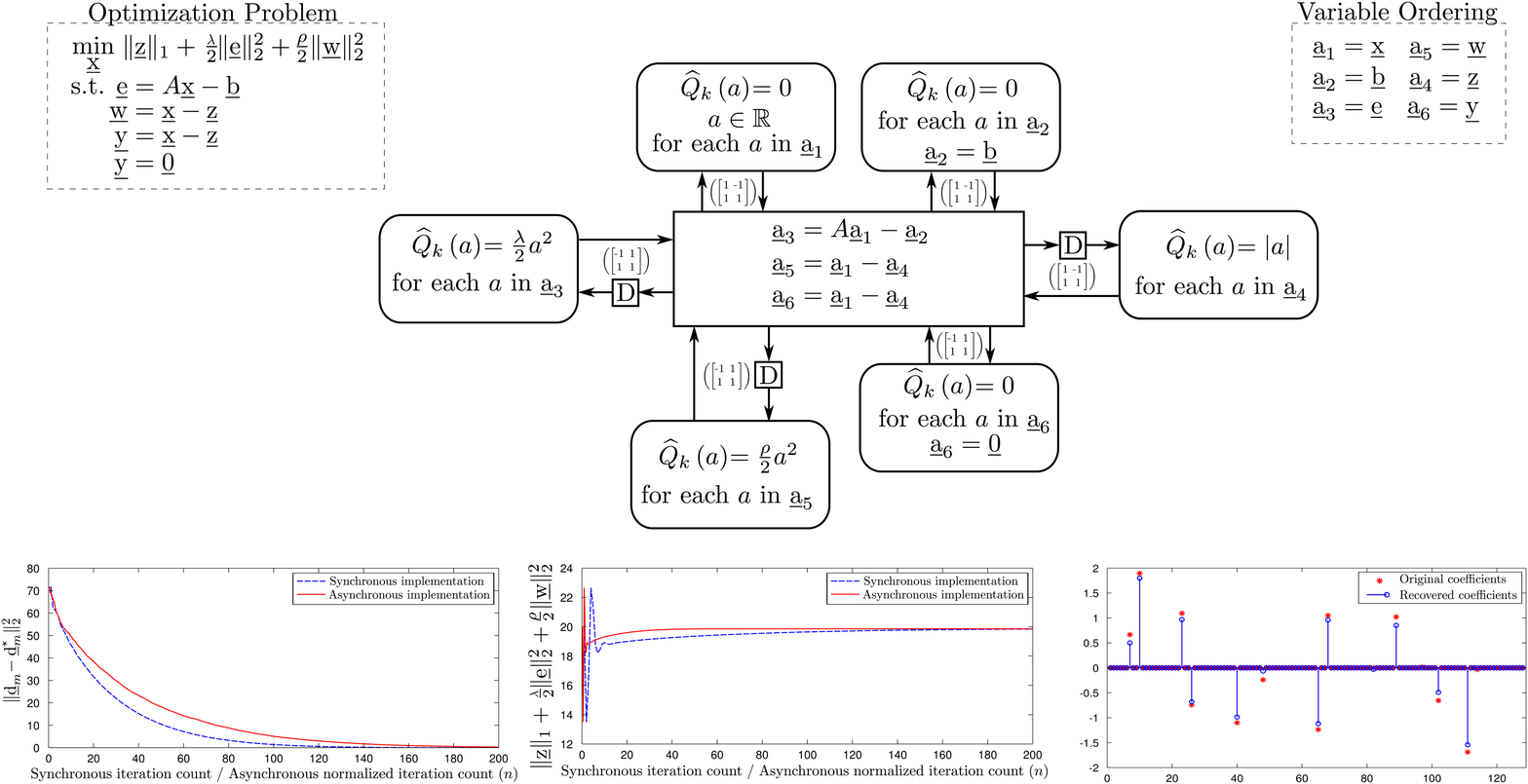,width=7in}
\caption{Signal processing architecture and numerical simulation for an augmented-cost LASSO problem, with the cost being augmented similarly to that of various ADMM formulations.\cite{BoydNealChu} The augmentation parameter is denoted $\rho$.  The parameter $\lambda$ is selected to trade off between the sparsity of the solution and the enforcement of the linear equality constraints. For the depicted solution $\lambda$ and $\rho$ are selected to be large.}\label{fig:lasso_augmented}
\vspace{-4mm}
\end{center}
\end{figure*}

\begin{figure*}[h!]
\begin{center}
\epsfig{file=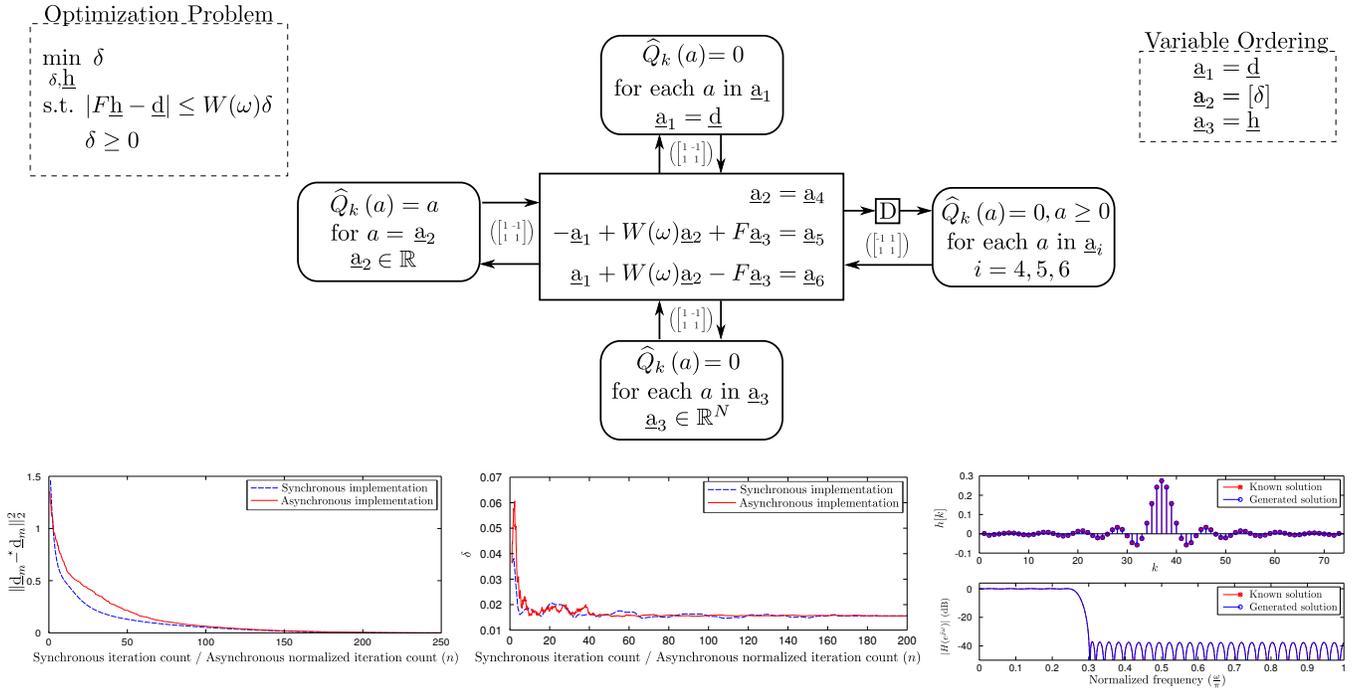,width=7in}
\caption{Signal processing architecture and numerical simulation corresponding to a minimax-optimal FIR filter design problem, specifically that of lowpass filter design. The obtained result is compared with a known solution from the Parks-McClellan algorithm.}\label{fig:minimax_filter_design}
\end{center}
\end{figure*}

\begin{figure*}[h!]
\begin{center}
\epsfig{file=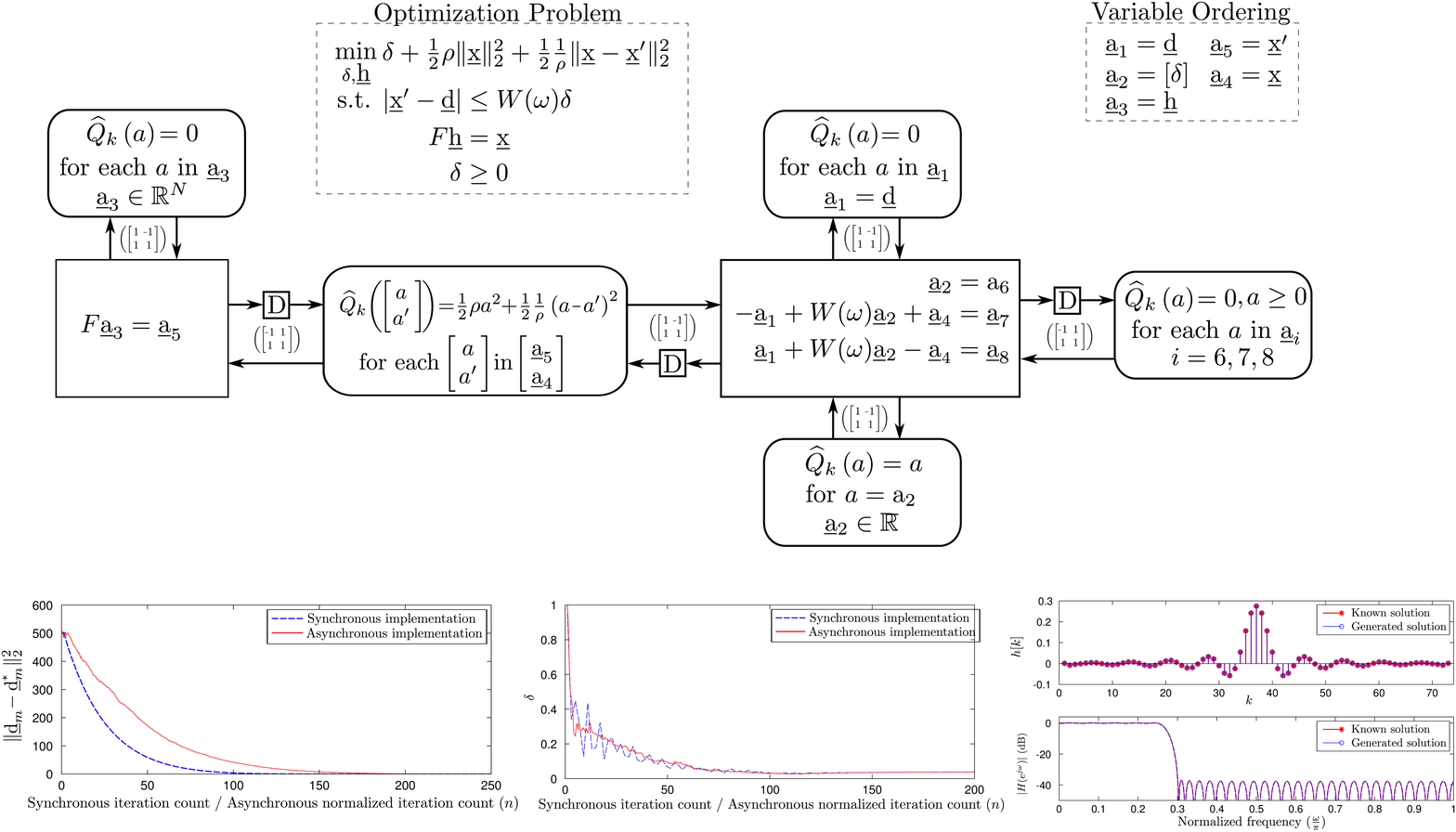,width=7in}
\caption{Alternative algorithm for minimax-optimal filter design, obtained by modification of the problem statement in Fig.~\ref{fig:minimax_filter_design} and intended to demonstrate that the presented framework can be used in designing a variety of distinct classes of algorithms. The parameter $\rho$ is selected to specify the relative enforcement of equality between the system variables loosely shared between the two linear interconnection elements. For the depicted solution $\rho$ is selected to be small, resulting in a very close approximation to the lowpass filter design problem in Fig.~\ref{fig:minimax_filter_design}.}\label{fig:minimax_filter_design_2}
\end{center}
\end{figure*}

\begin{figure*}[h!]
\begin{center}
\epsfig{file=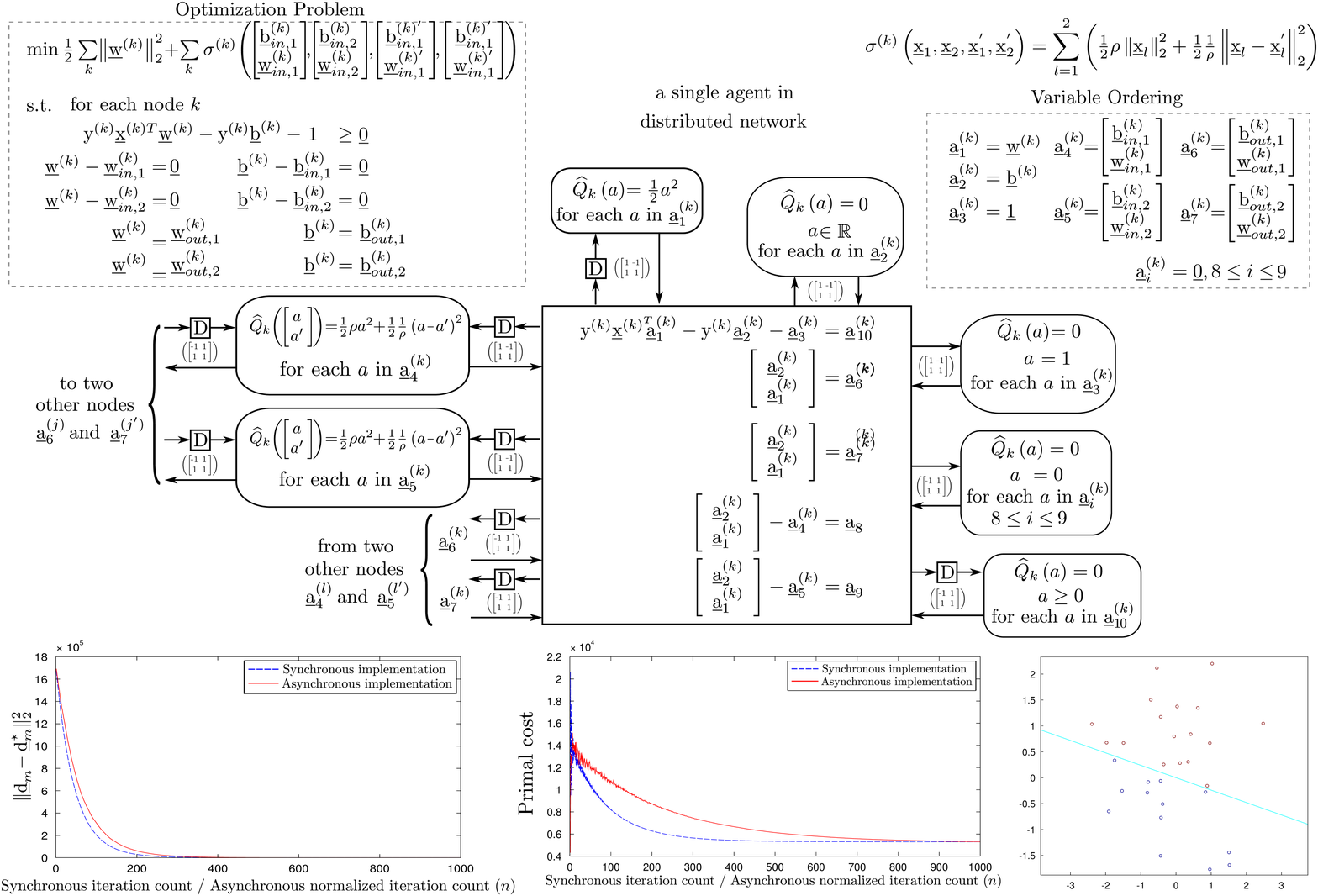,width=7in}
\caption{Signal processing architecture for a single agent in a connected graph implementing a decentralized algorithm for training a support vector machine classifier. The numerical simulation depicts a system involving $30$ such agents, each having knowledge of a single training vector. The parameter $\rho$ specifies the relative enforcement of equality for the system variables that are coupled between each agent in the graph. For the depicted solution $\rho$ is selected to be small, and the graph is known to be connected, with  each node as depicted above having exactly four incident connections.}\label{fig:svm}
\end{center}
\end{figure*}

\begin{figure*}[h!]
\begin{center}
\epsfig{file=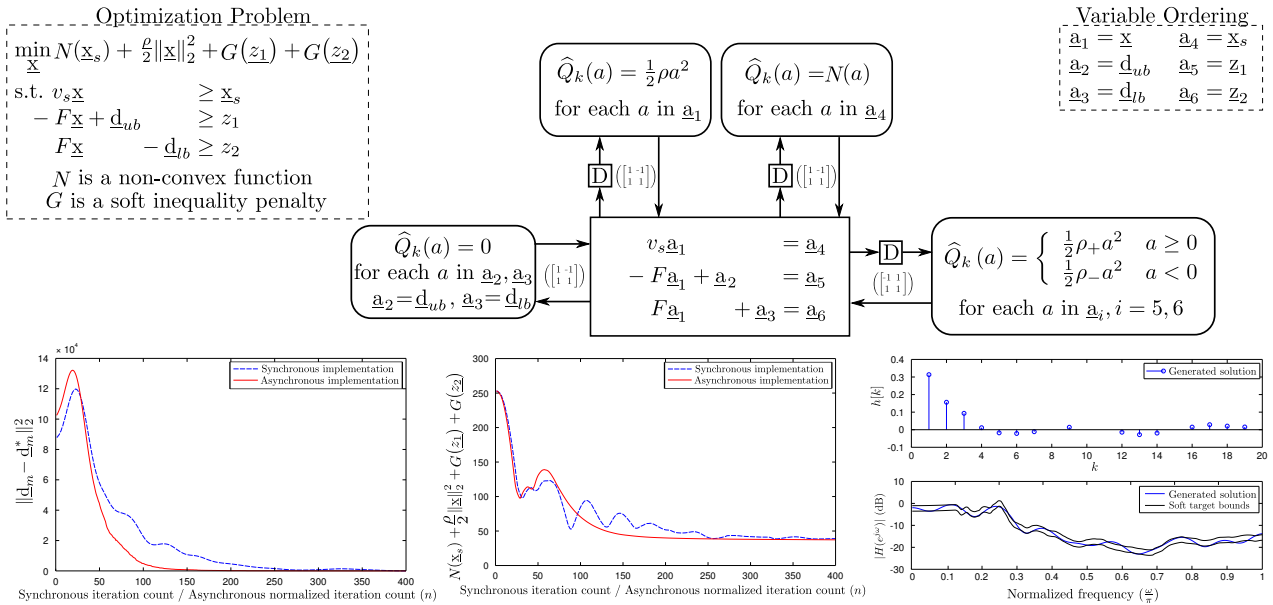,width=6.6in}
\caption{Signal processing architecture and numerical simulation corresponding to a nonconvex sparse filter design problem. The parameters $\rho$ and $v_{s}$ are respectively selected to specify the enforcement of the size of $\underline{\mbox{x}}$ and the width of the abrupt decrease in cost about $0$ for the nonconvex element.  $\rho_{+}$ and $\rho_{-}$ affect the enforcement of the soft inequality constraints. For the depicted solution $\rho$ and $\rho_{+}$ are selected to be small and $\rho_{-}$ and $v_{s}$ are selected to be large.}\label{fig:sparse_equalizer}
\end{center}
\end{figure*}

Fig.~\ref{fig:convergence}(a) summarizes the overall interconnection of elements composing the presented class of systems discussed in Part I \cite{BaranLahlouPartI}, with those maps $m_k(\cdot)$ corresponding to source relationships being written separately. Figs.~\ref{fig:convergence}(b)-(d) illustrate a set of manipulations useful in analyzing convergence, with Fig.~\ref{fig:convergence}(b) specifically depicting a solution to the transformed stationarity conditions.  The approach is to begin with the system in Fig.~\ref{fig:convergence}(a) and perform the additions and subtractions of $c^{\star}_i$ and $d^{\star}_i$ indicated in Fig.~\ref{fig:convergence}(c), obtaining Fig.~\ref{fig:convergence}(d) by identifying that Fig.~\ref{fig:convergence}(c) is a superposition of Figs.~\ref{fig:convergence}(b) and (d).

\vspace{-2mm}
There are various ways that the system in Fig.~\ref{fig:convergence}(d) can be used in determining sufficient conditions for convergence, a subset of which we outline here.  Generally, arguments for convergence utilizing Fig.~\ref{fig:convergence}(d) involve identifying conditions for which $||\underbar{d}_D||$ in this figure is strictly less than $||\underbar{d}'_m||$, except at $0$.  Using the definition of a source element in \cite{BaranLahlouPartI} and the fact that $G$ is a neutral map, i.e.~an orthonormal matrix, we conclude from Fig.~\ref{fig:convergence}(d) that
\begin{equation}
||\underbar{d}_D|| \leq ||\underbar{c}'_m||.
\end{equation}

If, for example, the solution to the transformed stationarity conditions $c^{\star}_i$ and $d^{\star}_i$ is known to be unique, and additionally if the collection of constitutive relations denoted $m(\cdot)$ is known to be dissipative about $\underbar{d}_m^{\star}$, then from Fig.~\ref{fig:convergence}(d) we conclude that $|| \underbar{c}'_m || < || \underbar{d}'_m ||$ except at $0$, resulting in
\begin{equation}
||\underbar{d}_D|| < || \underbar{d}'_m || \label{eq:normRedConv}
\end{equation}
except at $0$.  Eq.~\ref{eq:normRedConv} implies, for example, that coupling the constitutive relations denoted $m(\cdot)$ to the linear interconnection elements via deterministic vector delays, the discrete-time signal denoted $\underbar{d}'_m[n]$ will converge to $0$ and so the signal $\underbar{d}_m[n]$ will converge to $\underbar{d}^{\star}_m$.

The uniqueness of the stationarity conditions and the property of the constitutive relations being dissipative used in the preceding argument are not, however, strictly required.  A more general line of reasoning involves justifying Eq.~\ref{eq:normRedConv} in the vicinity of any such solution $c^{\star}_i$ and $d^{\star}_i$, for example by claiming that even if specific constitutive relations $m_k(\cdot)$ are norm-increasing, the overall interconnected system results in a map from $\underbar{d}'_m$ to $\underbar{d}_D$ that is norm-reducing in the vicinity of that solution.

Arguments for convergence involving essentially Eq.~\ref{eq:normRedConv} can also be applied in a straightforward way to systems utilizing asynchronous delays, modeled as discrete-time sample-and-hold systems triggered by independent Bernoulli processes.  In particular by taking the expected value of $|| \underbar{d}'_{(m)}[n] ||$, applying the law of total expectation, substituting in Eq.~\ref{eq:normRedConv}, and performing algebraic manipulations, it can be argued that $E[|| \underbar{d}'_{(m)}[n] ||]$ converges to $0$.  A more formal treatment of convergence is the subject of future work.

\FloatBarrier
\clearpage

\bibliographystyle{IEEEbib}
\nocite{*}
\bibliography{refs_part2}

\end{document}